Д. А. Пасечнюк[1]

[1]Президентский физико-математический лицей №239

# Планирование распределения ресурсов вышки мобильной связи

Рассматривается задача составления расписания распределения (временных) ресурсов базовой станции (сотовой вышки), осуществляющей взаимодействие клиентов (пользователей беспроводных мобильных устройств, имеющих доступ в интернет) и серверов, с которых они закачивают web-страницы (в общем случае файлы).

**Ключевые слова:** автоматическое составление расписания, планировщик, сотовая вышка, индексные стратегии, вычислительный эксперимент.

D. A. Pasechnyuk[1]

[1]Presidential Physics and Mathematics Lyceum №239

# Scheduling strategies for resource allocation in a cellular base station

The problem of scheduling the (time) resource allocation of a base station (cell tower) that interacts with clients (users of wireless mobile devices with Internet access) and servers from which they download web pages (files in general) is studied.

**Key words:** automated scheduling, scheduler, cellular base station, index strategies, computational experiment.

## 1. Введение

Задача составления расписания работы планировщика для вышки сотовой связи стала популярна в последнее десятилетие в связи с массовым использованием мобильных устройств для выхода в интернет и скачивания (просмотра) больших объемов информации. В литературе по данному направлению прочно укрепился подход, при котором работа планировщика задается некоторым индексом: в данный момент времени обслуживается тот клиент (пользователь мобильного устройства), у которого значения этого индекса наибольшее. Индексы легко рассчитываются по накопленной на вышке информации об истории обслуживания имеющихся сейчас клиентов и по известным текущим величинам пропускных способностей каналов передачи данных для каждого клиента. В литературе насчитываются десятки различных индексов, которые в зависимости от условий, оказываются оптимальными, т.е. приводящими к оптимальным в некотором смысле расписаниям. В действительности, строгих результатов в данной области мало. В основном ставятся численные эксперименты, позволяющие отбирать среди различных вариантов наилучшие.

В данной работе также рассматриваются только численные эксперименты. Прежде всего, в работе показано, что при естественных условиях с помощью (выпуклой) комбинации имеющихся индексов не удается получить лучший по свойствам индекс. Это довольно неожиданный результат. Начиная работу в этом направлении, автор был уверен в обратном. Также в работе были получены новые индексы, в некотором смысле обратные к общепринятым, которые показали хорошие результаты при использовании популярных у практиков критериев качества работы планировщика ALPT и logALPT.

Структура статьи следующая. В п. 2 приводится краткое описание технической стороны рассматриваемой задачи. В п. 3 задача ставится более формально, но все же, по-прежнему, не строго. Строго поставить задачу с критериями ALPT и logALPT, насколько нам известно, пока открытая задача. В п. 4 вводится понятие индексной стратегии планировщика и приводится небольшой обзор по наиболее конкурентным индексам, которые удалось найти в литературе. В п. 5 вводятся некоторые новые индексы, полученные с использованием информации о вероятностном распределении входных данных. Со всеми этими индексами проводились численные эксперименты на модельных данных, которые описаны в заключительном п. 6.

## 2. Краткое описание физики процесса

Вызвав в браузере web-страницу или просто начав скачивание какого-то файла с сотового устройства на базовую станцию отправляется соответствующий запрос, который впоследствии (посредством использующегося протокола, например, протокола TCP/IP) ретранслируется по сети и доходит (по указанному в запросе IP-адресу) на соответствующий сервер. С сервера начинается закачивание пакетов (отвечающих запросу) в буфер базовой станции. Буфер (ресурсы памяти), который базовая станция может выделить этим пакетам (поступающим с сервера), ограничен. Поэтому довольно типично закачивание больших файлов на сотовые устройства через базовую станцию осуществляется посредством нескольких заполнений/опустошений соответствующего буфера. Докачка с сервера продолжается только в момент, когда для этого освободится место (соответствующий буфер базовой станции полностью очистится), а это, в свою очередь, зависит от алгоритма, осуществляющего распределение ресурсов базовой станции. Характерное время, через которое начинается новая докачка в буфер базовой станции RTT, составляет 30 мс (миллисекунд). Есть известный закон, определяемый используемым протоколом TCP/IP, "геометрической прогрессии с насыщением", которому подчиняются размеры закачиваемых на базовую станцию порций файла. Базовая станция работает в определенном (ограниченном) диапазоне (несущих) частот. Этот диапазон разбит на интервалы частот, в каждом из которых базовая станция посылает (в каждую единицу времени $\sim 1$ мс) пакеты (данные) на одно из сотовых устройств. Таким образом, осуществляется постепенное опустошение буферов. Число пакетов (трафик) $r_k(t)$, которые можно в данный момент времени $t$ передать (послать) от базовой станции пользователю $k$ определяется его текущим местоположением (изменив это положение на несколько метров, пользователь заметно меняет интерференционную картину, складывающуюся в результате всевозможных отражений несущих волн от зданий и т.п.). Это приводит к тому, что $r_k(t)$ непостоянны во времени. Типично, что $r_k(t)$ колеблются вокруг некоторого среднего значения с амплитудой составляющей $\sim 30\%$ от этого среднего значения. Однако эти колебания носят существенно случайный (не регулярный) характер.

Все эти детали необходимо, конечно, учитывать при разработке алгоритма для планировщика, распределяющего ресурсы базовой станции среди клиентов. Для этого существуют специальные пакеты (симмуляторы) типа NS3. Однако в данной работе, как и в подавляющем большинстве других работ, разработка алгоритма для планировщика будет осуществляться при дополнительных упрощающих предположениях.

В частности, будем далее считать, что $\{r_k(t)\}_t$ является i.i.d. последовательностью (последовательностью независимых одинаково распределенных случайных величин). Величины $r_k(t)$ известны на базовой станции.

Аппаратура, установленная на базовой станции, вполне позволяет работать с поступающими данными (и хранить их) так, как это делает современный ноутбук. Доступным наблюдению на базовой станции являются объем уже переданного трафика на сотовое устройство ($R_k(t)$ — размер трафика, переданного к моменту $t$ пользователю $k$) и текущий размер буферов, отвечающих активным в каждый момент пользователям ($b_k(t)$ — размер буфера в момент $t$ пользователя $k$). Подчеркнем, что речь идет о размере буферов на базовой станции, но не на сервере. Размер буфера на сервере (объем того, что надо в итоге скачать), который может быть заметно больше, чем на базовой станции, к сожалению, не доступен наблюдению. Со временем ситуация может измениться, когда протокол передачи данных будет позволять передавать информацию о размере всего файла в первых пакетах и это можно будет считывать. Однако сейчас надо исходить из того, что базовая станция "не знает" сколько еще лежит на сервере данных, которые впоследствии надо будет закачать в соответствии с запросом пользователя. Объем файла, который хочет скачать пользователь $k$ обозначим через $A_k$. Для простоты мы отождествляем номер пользователя и номер файла. Понятно, что один пользователь может закачивать несколько файлов. Однако, не ограничивая общности, можно вторую сессию того же пользователя просто воспринимать как нового пользователя (если не пытаться учитывать в параметрической модели какие-то закономерности того, что это один и тот же человек – на данном этапе мы не будем на этом останавливаться).

Число клиентов (пользователей) доступных в данный момент также всегда известно и на каждого из них базовая станция заводит свой буфер. Число пользователей меняется со временем (кто-то приходит, кто-то уходит, когда его обслужили). Вся информация о появлении новых пользователей ($T_k^0$ — момент прихода пользователя $k$) и об уходе пользователей ($T_k^{end}$ — момент ухода пользователя $k$) на базовой станции известна.

## 3. Постановка задачи

Таким образом, на базовой станции хранится история $\{r_k(t)\}_{k,t}$, $\{b_k(t)\}_{k,t}$, $\{R_k(t)\}_{k,t}$, $\{T_k^0\}_k$, $\{T_k^{end}\}_k$, $\{A_k\}_k$, исходя из которой требуется построить оптимальное расписание распределения (временных) ресурсов базовой станции (будем обозначать это расписание буквой $x$) в смысле критерия ALPT (Application Level Perceived Throughput)

$$F(x) \stackrel{\text{п.н.}}{=} \lim_{T \to \infty} \frac{1}{N(T)} \sum_{k=1}^{N(T)} \frac{A_k}{T_k^{end}(x) - T_k^0} \to \max_x$$

или logALPT

$$F(x) \stackrel{\text{п.н.}}{=} \lim_{T \to \infty} \frac{1}{N(T)} \sum_{k=1}^{N(T)} \log\left(\frac{A_k}{T_k^{end}(x) - T_k^0}\right) \to \max_x.$$

В определении $F(x)$ предполагается "эргодичность", то есть существование для любого $x$ предела (почти наверное) в выражении в правой части равенства. К сожалению, строго здесь это определить не так просто, поскольку отсутствует полугрупповое свойство — функционал задачи не представляется как сумма значений некоторой функции на траектории некоторого случайного процесса в дискретном времени. Другими словами, классические результаты теории случайных процессов по эргодичности в данном случае не понятно как можно применить для обоснования существования такого предела. Однако далее наши рассуждения будут исходить из того, что функционал определен корректно.

Наша задача построить на основе доступной (на базовой станции) информации обучающийся алгоритм (с подкреплением), который распределяет ресурсы базовой станции оптимально, исходя из критерия максимума $F(x)$.

Детали постановки задачи можно найти в монографии [8].

В таком виде задача не сводится ни к чему классическому. В частности, даже к Markov Decision Process. Другими словами, в такой постановке невозможно даже записать уравнения Вальда–Беллмана, не говоря уже о том, чтобы его решить. И, тем более, отыскать решение в индексном виде (то есть вычислительно эффективно).

## 4. Краткий обзор литературы по близким постановкам (индексы)

Начиная с цикла работ Frank'a Kelly конца 90-х годов (основные работы [1, 2]) предполагается, что (конечное) число клиентов $N$ зафиксировано и не изменяется. Каждый клиент скачивает бесконечно большой файл ($A_k = \infty$, $k = 1, ..., N$). В качестве критерия качества используемого алгоритма составления расписания выбирается вогнутая функция $H(u)$ и рассматривается динамика throughput'ов $\theta_k(t) = R_k(t)/t$. При весьма общих условиях эта динамика будет эргодической (равномерно по используемым алгоритмам). Это означает, что в зависимости от выбранного алгоритма существует предел $u_k \stackrel{\text{п.н.}}{=} \lim_{t \to \infty} \theta_k(t)$. Перебирая всевозможные допустимые алгоритмы, мы породим выпуклое множество $U$. На этом множестве надо найти $u_* = \arg\max_{u \in U} H(u)$, и соответствующий ей алгоритм(-ы). Эта задача была решена в общем виде [2]. Однако здесь мы приведем решение в частном случае (наиболее популярном на практике): $H(u) = \sum_{k=1}^{N} \ln u_k$. В этом случае оптимальная (индексная) стратегия (см. определение индексной стратегии в [3]) имеет вид — выбирать для обслуживания в момент времени $t$ клиента с наибольшей величиной $r_k(t)/\theta_k(t)$. Стратегию

$$\arg\max_{k \in \left[\substack{\text{текущие} \\ \text{пользователи}}\right]} \frac{r_k(t)}{\theta_k(t-1)} \quad \textbf{(PFS)}$$

называют **Proportional-Fair Sharing** (PFS).

С технической точки зрения за такт времени базовая станция может обслуживать одновременно несколько клиентов. PFS (в варианте [2]) обобщается на этот случай. Однако, всегда можно путем дополнительного (условного) разделения (частотным образом) ресурсов внутри такта времени свести постановку к рассмотрению ситуации, когда за такт времени базовая станция обслуживает ровно одного клиента. Далее мы будем именно так и считать.

Если считать, что $r_k(t) \equiv r_k$ и в качестве критерия использовать

$$\tilde{F}(x) \stackrel{\text{п.н.}}{=} \lim_{T \to \infty} \frac{1}{T} \sum_{k=1}^{N(T)} \left( T_k^{end}(x) - T_k^0 \right) \to \min_x,$$

то функционал можно представить классическим (для динамического программирования) образом

$$\tilde{F}(x) \stackrel{\text{п.н.}}{=} \lim_{\beta \to 1-0} \frac{1}{1-\beta} \sum_{t=1}^{\infty} \sum_{k=1}^{\infty} \beta^t I \left[\substack{\text{k-й клиент появился,} \\ \text{но еще не полностью обслужен}}\right].$$

Оптимальный алгоритм здесь [7]

$$\arg\max_{k \in \left[\substack{\text{текущие} \\ \text{пользователи}}\right]} \frac{r_k(t)}{A_k - R_k(t-1)}.$$

К сожалению этот (индексный) алгоритм (его обычно называют **SRPT**) относится к классу opportunistic (size-aware, anticipating). Это означает, что алгоритм требует знания $A_k$ еще до момента полной загрузки клиентом файла. К сожалению, в реальности базовые станции узнают $A_k$ только в момент полной загрузки файла (и никак не раньше). С помощью индексов Гиттинса [3] можно предложить non-anticipating вариант. Слабым место здесь является не столько критерий [1], сколько предположение $r_k(t) \equiv r_k$, которое далеко от истины. Численные эксперименты с реальными данными показали, что такое предположение является слишком грубым.

В литературе можно найти такой вариант SRPT-индекса

$$\arg\max_{k \in \left[\substack{\text{текущие} \\ \text{пользователи}}\right]} \frac{r_k(t)}{b_k(t)} \quad \textbf{(SECTF)},$$

---

[1]Его (критерий/функционал) как раз можно доводить до интересующих нас вариантов, с точностью до усреднения обратных величин, путем различных взвешиваний — см., например, [4], [5] (и цитированную там в пп. 1, 7 литературу). Правда, при этом уже нельзя утверждать, что полученный метод будет оптимальным. Отметим, что в этих работах часто используются другие (более современные) индексы [6].

$$\arg\max_{k\in\left[\substack{\text{текущие}\\\text{пользователи}}\right]} \frac{r_k(t)}{R_k(t-1)} \quad \textbf{(DAS)}.$$

Имеются также другие индексы, которые в определенных ситуациях используют на практике. Например, Round Robin

$$\arg\max_{k\in\left[\substack{\text{текущие}\\\text{пользователи}}\right]} \frac{1}{t - T_k^0} \quad \textbf{(PS)},$$

$$\arg\max_{k\in\left[\substack{\text{текущие}\\\text{пользователи}}\right]} r_k(t) \quad \textbf{(MaxC/I)},$$

$$\arg\max_{k\in\left[\substack{\text{текущие}\\\text{пользователи}}\right]} \frac{r_k(t)}{t - T_k^0} \quad \textbf{(TAS)}.$$

Эти индексы имеют вполне понятный физический смысл (мы не будем здесь это обсуждать).

## 5. Предложенные индексы

Рассмотрим выражение для критерия ALPT

$$ALPT = \frac{1}{N}\sum_{k=1}^{N}\frac{A_k}{T_k^{end} - T_k^0},$$

после чего заменим выражение под знаком суммы на его приближенное значение

$$\frac{A_k}{T_k^{end} - T_k^0} = \frac{A_k}{(t - T_k^0) + (T_k^{end} - t)} \approx \frac{A_k}{(t - T_k^0) + \frac{A_k - R_k(t-1)}{\bar{r}_k}}. \tag{1}$$

В данном выражении фигурирует $A_k$, неизвестное нам до момента окончательной загрузки файла. Однако можно заменить $A_k$ в данной формуле на некоторое ожидаемое значение $A_k^{expected}$. Рассмотрим простой случай, когда $A_k$ есть случайная величина из распределения Парето $\Pi(A_k|m,\alpha)$, где $m$, $\alpha$ — некоторые известные константы. Тогда, по определению:

$$p(A_k) = \frac{\alpha m^\alpha}{A^{\alpha+1}} = \frac{const}{A^{\alpha+1}}.$$

Из этого выражения видим, что вероятность обратно пропорциональна $A^{\alpha+1}$, что и будет использовано в дальнейшем.

Известно, что к моменту времени $t$ на устройство уже была скачана некоторая часть файла размером $R_k(t-1)$. Следовательно, оценивая $A_k$ в момент времени $t$, следует предположить, что $\forall x < R_k(t-1)$ вероятность $\hat{p}(x) = 0$. То есть, в дальнейшем мы будем использовать некоторую вероятностную функцию $\hat{p}: \mathbb{R} \to [0,1]$, такую что $\forall x < R_k(t-1): \hat{p}(x) = 0$, а для $x \geq R_k(t-1)$ будет выполняться $\hat{p}(x) \sim \frac{1}{x^{\alpha+1}}$.

Пусть $C$ — коэффициент пропорциональности в выражении для $\hat{p}$. Потребуем, чтобы выполнялось

$$\int_{R_k(t-1)}^{\infty} p(A_k)dA_k = 1,$$

то есть

$$C\int_{R_k(t-1)}^{\infty}\frac{dA_k}{A^{\alpha+1}} = C\int_{R_k(t-1)}^{\infty}A^{-(\alpha+1)}dA_k = C\cdot\left(\frac{A_k^{-\alpha}}{-\alpha}\right)\bigg|_{R_k(t-1)}^{\infty} = \frac{C}{\alpha R_k^\alpha(t-1)} = 1,$$

откуда

$$C = \alpha R_k^\alpha(t-1).$$

Пусть

$$A_k^{expected} \stackrel{\text{def.}}{=} \int_{R_k(t-1)}^\infty A_k \hat{p}(A_k) dA_k.$$

Тогда, используя полученное выражение для $\hat{p}$ получаем:

$$A_k^{expected} = \int_{R_k(t-1)}^\infty \frac{\alpha R^\alpha}{A_k^\alpha} dA_k = \alpha R^\alpha \int_{R_k(t-1)}^\infty \frac{dA_k}{A^\alpha} = \alpha R^\alpha \left( \frac{A_k^{1-\alpha}}{1-\alpha} \right) \Bigg|_{R_k(t-1)}^\infty,$$

$$A_k^{expected} = \frac{\alpha}{\alpha-1} R_k^\alpha(t-1) \cdot R_k^{1-\alpha}(t-1) = \frac{\alpha}{\alpha-1} R_k(t-1).$$

Теперь подставим полученное ожидаемое значение $A_k^{expected}$ в оценку (1) выражения из формулы для ALPT, из чего получим некоторое ожидаемое значение критерия для $k$-го клиента.

$$\frac{\alpha}{\alpha-1} \cdot \frac{R_k(t-1)}{(t-T_k^0) + \frac{R_k(t-1)}{C_1 \bar{r}_k(t)}} \quad (\mathbf{T}),$$

которое может быть использовано в качестве индекса: максимизируя это выражение мы также будем максимизировать и одно из слагаемых критерия качества ALPT.

Глядя на выражение для T можно понять, что значение индекса будет увеличиваться в следующих случаях:

1) $R_k(t-1)$ возрастает;

2) $\bar{r}_k(t)$ возрастает;

3) $(t - T_k^0)$ убывает.

Тогда можно предложить упрощенный индекс, отражающий приведённые зависимости индекса T:

$$\frac{R_k(t-1)\bar{r}_k(t)}{t - T_k^0} \quad (\mathbf{TK}).$$

## 6. Вычислительный эксперимент

В рамках данной работы был поставлена серия вычислительных экспериментов, целью которых является сравнение различных подходов для построения индексных стратегий. Эксперименты проводились на синтетических данных и их задачей было ответить на следующие вопросы:

1) Насколько хорошо работают предложенные в данной работе индексы по сравнению с уже имеющимися?

2) Можно ли получить лучшее решение, если рассматривать не индексы по отдельности, а их линейные комбинации?

3) Можно ли получить лучшее решение, если рассматривать вероятностные комбинации индексов?

## 6.1. Описание данных

Для генерации данных был реализован эмулятор работы сотовой вышки и очереди клиентов. Код написан на языке Python 3.5 и находится в открытом доступе по ссылке[2]. В основе работы эмулятора лежит параметрическая модель генерации данных, базирующаяся на следующих предположениях:

- Времена между приходами двух подряд идущих клиентов $T_{k+1}^0 - T_k^0$ есть н. о. р. случайные величины из показательного распределения

$$p(x) = \text{Exp}(x|\lambda).$$

В наших экспериментах $\lambda = 0.09$, что соответствует средней нагрузке на систему в размере 11 клиентов в секунду.

- Размеры скачиваемых файлов $A_k$ есть н. о. р. случайные величины из смеси распределений Парето

$$p(x) = \sum_{i=1}^{k} p_i \Pi(x|m_i, \alpha),$$

где число элементов смеси соответствует количеству различных типов файлов, где каждый тип имеет свой характерный размер. В наших экспериментах использовалось 4 различных значения размеров файлов, соответствующих загрузке страницы с текстом, работе приложения, прослушивания аудиозаписи и просмотру видеоролика (в килобайтах): $m_1 = 500, m_2 = 5000, m_3 = 25000, m_4 = 62500$. Соответствующие значения вероятностей, $p_1 = 0.4, p_2 = 0.3, p_3 = 0.2, p_4 = 0.1$, параметр $\alpha = 5.5$.

- Пропускные способности устройств клиентов $r_k(t)$ есть н. о. р. случайные величины из равномерного распределения

$$p(x) = \mathbf{U}(x|a_k(t), b_k(t)).$$

В наших экспериментах каждому из пользователей в момент создания запроса к базовой станции присваевается средняя пропускная способность $\bar{r}_k$, являющаяся случайной величиной из равномерного распределения

$$p(x) = \mathbf{U}\left(x \,\middle|\, \frac{\lambda \bar{A}}{3}, 3\lambda \bar{A}\right).$$

И для каждого определенного пользователя

$$a_k(t) = 0.7 \cdot \frac{3}{2} \left(\sin\left(10^{-4} \cdot 5 + 0.1\right) + 1\right) \cdot \bar{r}_k,$$

$$b_k(t) = 1.3 \cdot \frac{3}{2} \left(\sin\left(10^{-4} \cdot 5 + 0.1\right) + 1\right) \cdot \bar{r}_k.$$

## 6.2. Сравнение различных индексов

В первом эксперименте мы оцениваем качество работы планировщиков, построенных на различных индексах, и сравниваем как качество работы предложенных в работе индексов соотносится с ними. Максимизируемым функционалом качества был выбран logALPT (аналогичные результаты наблюдались и для ALPT):

$$\text{logALPT} = \frac{1}{N} \sum_{k=1}^{N} \log\left(\frac{A_k}{T_k^{\text{end}} - T_k^0}\right).$$

---

[2]GitHub: https://github.com/dmivilensky/Throughput-maximization-for-flows-with-fine-input-structure

В табл. 1.1 приведены значения критерия logALPT для различных индексов в эксперименте на синтетических данных. Как видно, предложенные индексы T и TK реализуют лучшие стратегии планировщика, чем все остальные индексы. В формуле, задающей индекс T, значение константы $C = \frac{0.6}{\ln(13/7)}$.

| Индекс | Формула | logALPT |
|---|---|---|
| **T** | $\frac{R_k(t-1)}{t-T_k^0 + \frac{R_k(t-1)}{C\bar{r}_k(t)}}$ | **4.01 ± 0.09** |
| **TK** | $\frac{r_k(t)R_k(t-1)}{t-T_k^0}$ | **3.93 ± 0.03** |
| RoundRobin | $\frac{1}{t-T_k^0}$ | 3.69 ± 0.03 |
| TAS | $\frac{r_k(t)}{t-T_k^0}$ | 3.50 ± 0.05 |
| Max C/I | $r_k(t)$ | 1.38 ± 0.24 |
| DAS | $\frac{r_k(t)}{R_k(t-1)}$ | 0.71 ± 0.14 |
| PF | $\frac{r_k(t)(t-T_k^0)}{R_k(t-1)}$ | −1.16 ± 0.05 |

Таблица 1.1. Значение функционала logALPT для различных индексных стратегий.

### 6.3. Линейная комбинация индексов

В данном эксперименте новый индекс задаётся как линейная коминация нескольких базовых индексов. Т.к. каждый индекс представляет собой некую оценку приоритета обслуживания клиента, их можно объединить в взвешенную линейную комбинацию, которая будет учитывать особенности сразу нескольких индексов:

$$I_{\text{lin}}(I_1, I_2, \ldots, I_k) = \sum_{i=1}^{k} w_i I_i, \quad w_1, \ldots, w_k \geq 0.$$

На рис. 1 показан график зависимости критерия качества logALPT от весов линейной комбинации (в данном случае это один параметр $\alpha$, т.к. индексы, отличающиеся в константное число раз дают одну и ту же стратегию планировщика и мы можем положить вес первого индекса равным 1). Лучшее значение функционала достигается при $\alpha = 0$, что соответствует индексу $I_{\text{lin}} = I_{\text{TAS}}$. Линейная комбинация индексов не даёт никакого выигрыша.

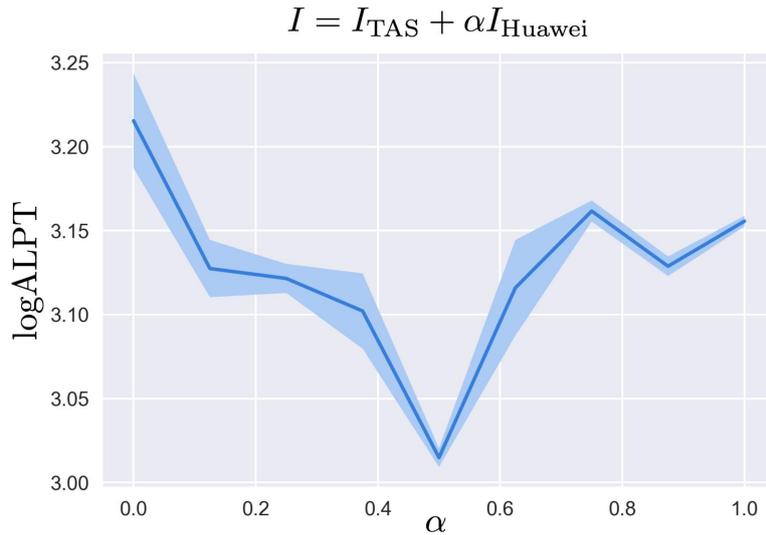

Рис. 1. Значение критерия logALPT для линейной комбинации двух индексов: TAS и DAS. Как видно, наибольшее значение достигается на границе, что соответствует использованию одного лишь индекса TAS.

### 6.4. Вероятностная комбинация индексов

В данном эксперименте индексная стратегия представляет собой вероятностную смесь нескольких уже известных индексов. Каждый раз перед тем как сделать выбор какого клиента обслужить, разыгрывается случайная величина, определяющая индекс, который будет использоваться на данном шаге. Получаемая новая индексная стратегия может быть формально записана как:

$$I_{\text{prob}}(I_1, I_2, \ldots, I_k) = \begin{cases} I_1, \text{ с вероятностью } p_1; \\ I_2, \text{ с вероятностью } p_2; \\ \ldots \\ I_k, \text{ с вероятностью } p_k. \end{cases} \quad \sum_{i=1}^{k} p_i = 1.$$

На рис. 2 показан график зависимости функционала logALPT от вероятностей элементов смеси $p_1, p_2, p_3$ для комбинации из трёх индексов: T, TAS и DAS. Лучшее значение функционала достигается тогда, когда мы на каждом шаге эмуляции выбираем индекс T и никогда не выбираем TAS или DAS. Вероятностная комбинация индексов не позволяет улучшить качество работы планировщика.

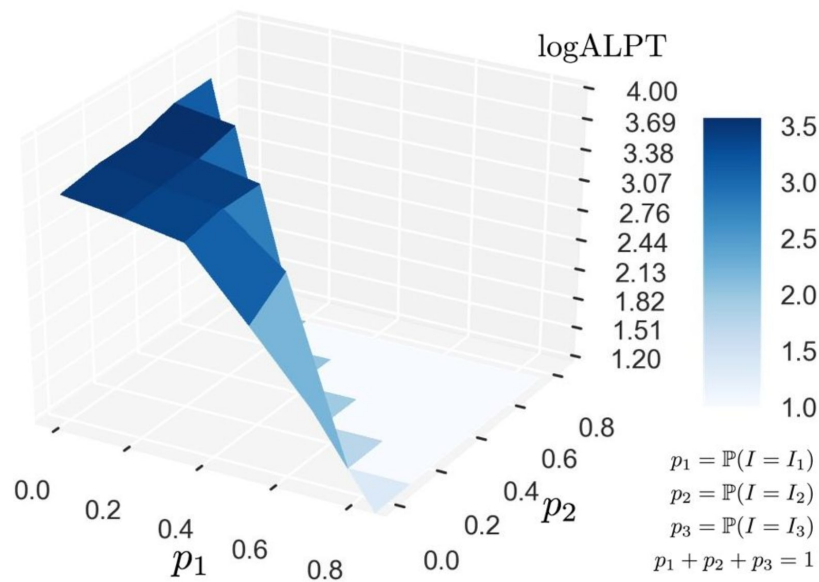

Рис. 2. Значение критерия logALPT для вероятностной комбинации трёх индексов: T, TAS и DAS. Как видно, наибольшее значение достигается на границе, что соответствует постоянному выбору одного индекса (T) с вероятностью единица.

## 7. Заключение

В работе исследуются различные индексные стратегии составления расписания обслуживания клиентов базовой станции (вышки сотовой связи). Показано, что при естественных условиях с помощью (выпуклой) комбинации индексов не удается получить лучший по свойствам индекс. Также в работе были предложены новые индексы, в некотором смысле обратные к общепринятым, которые показали хорошие результаты при использовании популярных у практиков критериев качества составленного расписания: ALPT и logALPT.

## 8. Благодарности

## Литература